\documentclass[letterpaper, 10pt, conference]{IEEEconf}
\IEEEoverridecommandlockouts 

\usepackage{amsmath,amssymb}
\usepackage{epsfig,subfigure,url,warmread}
\usepackage[all,import]{xy}

\newcommand{\norm}[1]{\ensuremath{\left\| #1 \right\|}}
\newcommand{\bracket}[1]{\ensuremath{\left[ #1 \right]}}
\newcommand{\braces}[1]{\ensuremath{\left\{ #1 \right\}}}
\newcommand{\parenth}[1]{\ensuremath{\left( #1 \right)}}
\newcommand{\refeqn}[1]{(\ref{eqn:#1})}
\newcommand{\reffig}[1]{Fig. \ref{fig:#1}}
\newcommand{\tr}[1]{\mbox{tr}\ensuremath{\negthickspace\bracket{#1}}}
\newcommand{\deriv}[2]{\ensuremath{\frac{\partial #1}{\partial #2}}}
\newcommand{\SO}{\ensuremath{\mathrm{SO(3)}}}
\newcommand{\T}{\ensuremath{\mathrm{T}}}
\newcommand{\so}{\ensuremath{\mathfrak{so}(3)}}

\renewcommand{\Re}{\ensuremath{\mathbb{R}}}
\renewcommand{\S}{\ensuremath{\mathbb{S}}}

\renewcommand{\theenumi}{\roman{enumi}}

\title{\LARGE \bf
Deterministic Global Attitude Estimation}

\author{Taeyoung Lee\authorrefmark{1}\authorrefmark{2}, Amit Sanyal, Melvin Leok\authorrefmark{1}, and N. Harris McClamroch\authorrefmark{2}%
\thanks{Taeyoung Lee, Aerospace Engineering, University of Michigan, Ann Arbor, MI 48109 {\tt tylee@umich.edu}}%
\thanks{Amit Sanyal, Mechanical and Aerospace Engineering, Arizona State University, Tempe, AZ
85287 {\tt sanyal@asu.edu}}%
\thanks{Melvin Leok, Mathematics, Purdue University, West Lafayette, IN 47907 {\tt mleok@math.purdue.edu}}%
\thanks{N. Harris McClamroch, Aerospace Engineering, University of Michigan, Ann Arbor, MI 48109 {\tt nhm@umich.edu}}%
\thanks{\textsuperscript{\footnotesize\ensuremath{*}}This research has been supported in part by NSF under grant DMS-0504747, and by a grant from the Rackham Graduate School, University of Michigan.}
\thanks{\textsuperscript{\footnotesize\ensuremath{\dagger}}This research has been supported in part by NSF under grant ECS-0244977.}
}

\begin{document}
\maketitle \thispagestyle{empty} \pagestyle{empty}

\begin{abstract}
A deterministic attitude estimation problem for a rigid body in an
attitude dependent potential field with bounded measurement errors
is studied. An attitude estimation scheme that does not use generalized
coordinate representations of the attitude is presented here. Assuming
that the initial attitude, angular velocity and measurement noise lie
within given ellipsoidal bounds, an uncertainty ellipsoid that bounds
the attitude and the angular velocity of the rigid body is obtained.
The center of the uncertainty ellipsoid provides point estimates, and
its size gives the accuracy of the estimates. The point estimates and
the uncertainty ellipsoids are propagated using a Lie group variational
integrator and its linearization, respectively. The estimation scheme
is optimal in the sense that the attitude estimation error and the size
of the uncertainty ellipsoid is minimized at each measurement instant,
and it is global since the attitude is represented by a rotation matrix.
\end{abstract}

\section{Introduction}
Attitude estimation is often a prerequisite for controlling
aerospace and underwater vehicles, mobile robots, and other
mechanical systems moving in space. The attitude determination
problem for a rigid body from vector measurements was first posed in
\cite{jo:wahba}. A sample of the literature in attitude estimation
can be found in \cite{BaOs.AES85, Shu.JAS90, ReHu.Au04}.

Most existing attitude estimation schemes use generalized coordinates
to represent the attitude. As is well known, minimal coordinate
representations of the rotation group, like Euler angles, Rodrigues
parameters, and modified Rodrigues parameters, lead to geometric or
kinematic singularities. Non-minimal coordinate representations, like
quaternions used in the quaternion estimation (QUEST) algorithm and
its several variants (\cite{Shu.JAS90,Pi.JGCD00}), have their own
associated problems. Besides the extra unit norm constraint one
needs to impose on the quaternion, the quaternion representation,
which is diffeomorphic to $\mathrm{SU(2)}$, double covers $\SO$. As such,
it has an inevitable ambiguity in expressing the attitude.

A stochastic state estimator requires probabilistic models for the
state uncertainty and the noise. However, statistical properties of
the uncertainty and the noise are often not available. We usually
make statistical assumptions on disturbance and noise in order to
make the estimation problem mathematically tractable. In many
practical situations such idealized assumptions are not appropriate,
and this may cause poor estimation
performance~\cite{jo:TheSkaSou.IEEETSP94}.

An alternative deterministic approach is to specify bounds on the
uncertainty and the measurement noise without an assumption on their
distribution. Noise bounds are available in many cases, and
deterministic estimation is robust to the noise distribution. An
efficient but flexible way to describe the bounds is using
ellipsoidal sets, referred to as uncertainty ellipsoids. The idea of
the deterministic estimation process is based on set theory results
developed in~\cite{jo:Sc1968}; optimal deterministic estimation
problems are studied in~\cite{jo:MaNo1996} and \cite{jo:DuWaPo2001}
using uncertainty ellipsoids.

In this paper, we study attitude estimation problems for the
uncontrolled dynamics of a rigid body in an attitude-dependent
potential field using uncertainty ellipsoids. The estimation scheme
we present has the following important features: (1) the attitude is
globally represented by a rotation matrix without using coordinates,
(2) the deterministic estimator is distinguished from a Kalman or
extended Kalman filter, (3) the measurement errors are assumed to be
bounded but there is no restriction on their distribution, and (4)
the estimates are optimal in the sense that the size of uncertainty
is minimized at each estimation step.

This paper is organized as follows. The attitude determination
problem from vector observations is introduced in Section II. The
attitude estimation problem is formulated in Section III, and the
attitude estimation scheme with angular velocity measurements is
developed in Section IV. Numerical examples are presented in Section
V.

\section{Attitude Determination from vector observations}

Attitude of a rigid body is defined as the orientation of a body
fixed frame with respect to a reference frame. It is represented by
a rotation matrix that is a $3\times 3$ orthogonal matrix with
determinant 1. Rotation matrices have a group structure denoted by
$\SO$. The group action of $\SO$ on $\Re^3$ transforms a vector
represented in the body frame into the reference frame.
In the attitude estimation problem, we measure directions in the
body frame to fixed points with known directions in the reference
frame. The directions in the body frame are transformed into the
known reference directions by pre-multiplying by the rotation
matrix defining the attitude of the rigid body. The rotation matrix
can be estimated by minimizing an error between the transformed
measured directions and the known reference directions.

We denote the $i$th known direction vector in the reference frame
as $e^i\in\S^2$, and the corresponding vector represented in the
body frame as $b^i\in\S^2$. These direction vectors are normalized
to have unit lengths. The $e^i$ and $b^i$ vectors are related by a
rotation matrix $R\in\SO$ that defines the attitude of the rigid
body;
$e^i = R b^i$,
for all $i\in\braces{1,2,\cdots,m}$, where $m$ is the number of
measurements.

We assume that $b^i$ is measured by sensors in the body frame. Let
the measured direction vector be $\tilde{b}^i\in\S^2$, which contains
sensor errors, and denote the estimated rotation matrix by $\hat{R}
\in\SO$. The estimation error is given by $e^i-\hat{R}\tilde{b}^i$.
The attitude determination problem consists of finding $\hat{R}\in
\SO$ such that the weighted $2$ $\mathrm{norm}$ of those errors is
minimized.
\begin{align}
\min_{\hat{R}} \mathcal J & =\frac{1}{2}\sum_{i=1}^m w_i
(e^i-\hat{R}\tilde b^i)^T (e^i-\hat{R}\tilde b^i),\label{eqn:wahba}\\
&\text{subject to } \hat{R}\in\SO,\nonumber
\end{align}
where $E=\bracket{e^1,\cdots,e^m}\in\Re^{3\times m}$, $\tilde
B=\bracket{\tilde b^1,\cdots,\tilde b^m}\in\Re^{3\times m}$, and
$W=\mathrm{diag}\bracket{w^1,\cdots,w^m}\in\Re^{m\times m}$ is a
weighting factor for each measurement.

This problem is known as Wahba's problem~\cite{jo:wahba}. The
original solution of Wahba's problem is given in~\cite{jo:solwahba},
and a solution expressed in terms of quaternions (QUEST) is
presented in~\cite{jo:ShOh1981}. We use the solution expressed in
terms of a rotation matrix without using generalized
coordinates~\cite{pro:san2006}. A necessary and sufficient condition
for optimality of \refeqn{wahba} is given by
\begin{align}
\hat{R}=S L \in\SO,\;\ S=S^T >0 \label{eqn:san2006}
\end{align}
where $L=EW\tilde B^T \in\Re^{3\times 3}$ is non-singular. The unique
solution of \refeqn{san2006} is obtained by QR factorization of
$L=Q_qQ_r$
\begin{gather}
\hat{R}=\parenth{Q_q\sqrt{(Q_rQ_r^T)^{-1}}Q_q^T} L,\label{eqn:Rhat}
\end{gather}
where $Q_q\in\SO$, $Q_r\in\Re^{3\times 3}$ is an upper triangular
matrix, and the symmetric positive definite (principal) square root
is used. Equation \refeqn{Rhat} is the unique solution of Wahba's
problem~\cite{pro:san2006}.

\section{Attitude Estimation Problem formulation}

\subsection{State bounding estimation}

We use deterministic state bounding estimation using ellipsoidal sets,
referred to as uncertainty ellipsoids, to describe state uncertainty
and measurement noise. 
This deterministic estimation procedure has steps similar to those in
the Kalman filter, and is illustrated in \reffig{ue}. The left figure
shows time evolution of an uncertainty ellipsoid,
and the right figure shows a cross section at a fixed measurement
instant. At the $k$th time step, the state is bounded by an uncertainty
ellipsoid centered at $\hat{x}_k$. This initial ellipsoid is propagated
through time. Suppose that the state is measured next at the $(k+l)$th
time step, when the predicted uncertainty ellipsoid is centered at
$\hat{x}_{k+l}^f$. At this instant, the measurement uncertainty ellipsoid
is centered at $\hat{x}_{k+l}^m$. The actual state then lies in the
intersection of the two ellipsoids. In the estimation process, we find
a new ellipsoid that contains this intersection, as shown in the right
figure. The center of the new ellipsoid, $\hat{x}_{k+l}$ gives a point
estimate of the state at time step $k+l$, and the magnitude of the new
uncertainty ellipsoid measures the estimation accuracy.
The deterministic estimates are optimal in the sense that the sizes
of the ellipsoids are minimized.

\renewcommand{\xyWARMinclude}[1]{\includegraphics[width=0.30\textwidth]{#1}}
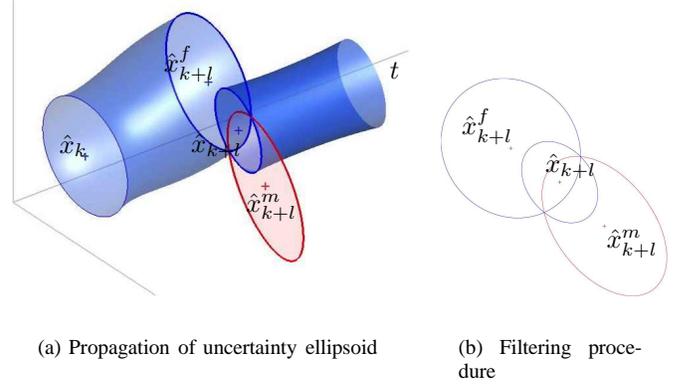
\begin{figure}[t]
    \centerline{\subfigure[Propagation of uncertainty ellipsoid]{
    $$\begin{xy}
    \xyWARMprocessEPS{tube}{eps}
    \xyMarkedImport{}
    \xyMarkedMathPoints{1-5}
    \end{xy}$$}
    \renewcommand{\xyWARMinclude}[1]{\includegraphics[width=0.17\textwidth]{#1}}
    \subfigure[Filtering procedure]{
        $$\begin{xy}
    \xyWARMprocessEPS{tubesection}{eps}
    \xyMarkedImport{}
    \xyMarkedMathPoints{1-3}
    \end{xy}$$}}
    \caption{Uncertainty ellipsoids}\label{fig:ue}
\end{figure}

\subsection{Equations of motion}
We consider estimation of the attitude dynamics of a rigid body in
the presence of an attitude dependent potential,
$U(\cdot):\SO\mapsto\Re$, $R\in\SO$. Systems that can be so modeled
include a free rigid body, spacecraft on a circular orbit with
gravity gradient effects~\cite{pro:acc06}, or a 3D
pendulum~\cite{pro:cca05}. The continuous equations of motion are
\begin{gather}
J\dot\Omega + \Omega\times J\Omega = M,\\
\dot{R} = R S(\Omega),\label{eqn:Rdot}
\end{gather}
where $J\in\Re^{3\times 3}$ is the moment of inertia matrix of the
rigid body, $\Omega\in\Re^3$ is the angular velocity of the body
expressed in the body fixed frame, and $S(\cdot):\Re^3\mapsto \so$
is a skew mapping defined by $S(x)y=x\times y$ for all $x,y\in\Re^3$.
The vector $M\in\Re^3$ is the moment due to the potential,
determined by $S(M)=\deriv{U}{R}^TR-R^T\deriv{U}{R}$,
or more explicitly,
\begin{gather}
M=r_1\times v_{r_1} + r_2\times v_{r_2} +r_3\times v_{r_3},
\end{gather}
where $r_i,v_{r_i}\in\Re^{1\times 3}$ are the $i$th row vectors of
$R$ and $\deriv{U}{R}$, respectively.

General numerical integration methods like the popular Runge-Kutta
schemes, typically preserve neither first integrals nor the
characteristics of the configuration space, $\SO$. In particular,
the orthogonal structure of the rotation matrices is not
preserved numerically. It is often proposed to parameterize
\refeqn{Rdot} by Euler angles or quaternions instead of integrating
\refeqn{Rdot} directly. However, Euler angles yield only local
representations of the attitude and they have singularities. Unit
quaternions do not exhibit singularities, but they have the manifold
structure of the three sphere $\S^3$, and double cover $\SO$.
Consequently, the unit quaternion representing the attitude is
inevitably ambiguous. In addition, general numerical integration
methods do not preserve the unit length constraint. Therefore,
quaternions have the same numerical drift problem as rotation
matrices.

Lie group variational integrators preserve the group structure
without the use of local charts, reprojection, or constraints,
they are symplectic and momentum preserving, and they exhibit good
energy behavior for an exponentially long time period. The following
Lie group variational integrator for the attitude dynamics of a
rigid body is presented in~\cite{pro:cca05}:
\begin{gather}
h S(J\Omega_k+\frac{h}{2} M_k) = F_k J_d - J_dF_k^T,\label{eqn:findf0}\\
R_{k+1} = R_k F_k,\label{eqn:updateR0}\\
J\Omega_{k+1} = F_k^T J\Omega_k +\frac{h}{2} F_k^T M_k
+\frac{h}{2}M_{k+1},\label{eqn:updatew0}
\end{gather}
where $J_d\in\Re^{3\times 3}$ is a nonstandard moment of inertia matrix
defined by $J_d=\frac{1}{2}\mathrm{tr}\!\bracket{J}I_{3\times 3}-J$,
and $F_k\in\SO$ is the relative attitude over an integration step.
The constant $h\in\Re$ is the integration step size, and the
subscript $k$ denotes the $k$th integration step. This integrator
yields a map $(R_k,\Omega_k)\mapsto(R_{k+1},\Omega_{k+1})$ by
solving \refeqn{findf0} to obtain $F_k\in\SO$ and substituting it
into \refeqn{updateR0} and \refeqn{updatew0} to obtain $R_{k+1}$ and
$\Omega_{k+1}$.

It preserves the orthogonal structure of $\SO$ because the rotation
matrix is updated by a product of two rotation matrices in
\refeqn{updateR0}. Since this integrator is obtained from a discrete
variational principle, it is symplectic, momentum preserving, and has good energy behavior, properties that are characteristic
of variational integrators. 

\subsection{Uncertainty Ellipsoid}
An uncertainty ellipsoid in $\Re^n$ is defined as
\begin{align}
    \mathcal{E}_{\Re^n}(\hat x,P)=\braces{x\in\Re^n \,\Big|\,
    (x-\hat{x})^T P^{-1}(x-\hat{x})\leq
    1},
\end{align}
where $\hat{x}\in\Re^n$, and $P\in\Re^{n\times n}$ is a symmetric
positive definite matrix. We call $\hat x$ the center of the
uncertainty ellipsoid, and $P$ is the uncertainty matrix that
determines the size and the shape of the uncertainty ellipsoid. The
size of an uncertainty ellipsoid is measured by $\mathrm{tr}\!\bracket{P}$
which is the sum of the squares of the semi principal axes of the
ellipsoid.

The state evolves in the 6 dimensional tangent bundle, $\T\SO$. We
identify $\T\SO$ with $\SO\times\so$ by left trivialization, and we
identify $\so$ with $\Re^3$ by the isomorphism $S(\cdot)$. The
uncertainty ellipsoid centered at $(\hat{R},\hat\Omega)\in\T\SO$ is
induced from an uncertainty ellipsoid in $\Re^6$;
\begin{align}
    \mathcal{E}(\hat{R},\hat{\Omega},P)
        & = \braces{R\in\SO,\,\Omega\in\Re^3 \,\Big|\,
    \begin{bmatrix}\zeta\\\delta\Omega\end{bmatrix}\in\mathcal{E}_{\Re^6}(0_{6},P)},\label{eqn:ueso}
\end{align}
where $S(\zeta)=\mathrm{logm} \parenth{\hat{R}^T R}\in\so$,
$\delta\Omega=\Omega-\hat{\Omega}\in\Re^3$, and $P\in\Re^{6\times
6}$ is a symmetric positive definite matrix. An element $(R,\Omega)
\in\mathcal{E}(\hat{R},\hat{\Omega},P)$ can be written as
\begin{align*}
    R = \hat{R} e^{ S(\zeta)},\quad
    \Omega = \hat{\Omega} + \delta \Omega,
\end{align*}
for some $x=\bracket{\zeta ;\,\delta\Omega}\in\Re^6$ satisfying
$x^TP^{-1}x\leq 1$.

\subsection{Uncertainty model}
We define the measurement error models for the direction vector and
for the angular velocity.
The measurement error is modeled by rotation of the measured
direction;
\begin{align}
{b}^i& = e^{S(\nu^i)} \tilde b^i,\nonumber\\
& \simeq \tilde b^i + S(\nu^i)\tilde b^i,\label{eqn:bi}
\end{align}
where $\nu^i\in\Re^3$ is the sensor error, which represents the
Euler axis of rotation vector from $\tilde b^i$ to $b^i$, and
$\norm{\nu^i}$ is the corresponding rotation angle in radians. The approximation is obtained by assuming that the measurement error is small.

The angular velocity measurement errors are modeled as
\begin{align}
\Omega_k=\tilde\Omega_k + \upsilon_k,\label{eqn:Omega}
\end{align}
where $\tilde\Omega_k \in\Re^3$ is the measured angular velocity,
and $\upsilon_k\in\Re^3$ is an additive error.

We assume that the initial conditions and the sensor noise are
bounded by prescribed uncertainty ellipsoids.
\begin{gather}
    (R_0,\Omega_0)\in\mathcal{E}(\hat{R}_0,\hat{\Omega}_0,P_0),\label{eqn:P0}\\
    \nu_k^i\in\mathcal{E}_{\Re^3}(0,S^i_k),\label{eqn:Sk}\\
    \upsilon_k\in\mathcal{E}_{\Re^3}(0,T_k)\label{eqn:Tk},
\end{gather}
where $P_0\in\Re^{6\times 6}$, $S_k^i, T_k \in\Re^{3\times 3}$ are
symmetric positive definite matrices that define the shape and the
size of the uncertainty ellipsoids.

\section{Attitude Estimation with Attitude and Angular Velocity Measurements}
In this section, we develop a deterministic estimator for the
attitude and the angular velocity of a rigid body assuming that both
attitude and angular velocity measurements are available. The estimator
consists of three stages; flow update, measurement update, and filtering.
The flow update predicts the uncertainty ellipsoid in the future. The
measurement update obtains an uncertainty ellipsoid using new measurements
and the sensor error model. Filtering obtains a new uncertainty
ellipsoid compatible with the predicted and the measured uncertainty
ellipsoids.

The subscript $k$ denotes the $k$th discrete index, and the
superscript $i$ denotes $i$th directional sensor. The superscripts
$f$ and $m$ denote the variables related to the flow update and the
measurement update, respectively. $\tilde\cdot$ denotes a measured
variable, and $\hat\cdot$ denotes an estimated variable.

\subsection{Flow update}
Suppose that the attitude and the angular momentum at the $k$th step
lie in a given uncertainty ellipsoid
\begin{align*}
    (R_k,\Omega_k)\in\mathcal{E}(\hat{R}_k,\hat{\Omega}_k,P_k),
\end{align*}
and that new measurements are taken at $(k+l)$th time step.
Flow update predicts the center and the uncertainty matrix that
define the uncertainty ellipsoid at the $(k+l)$th step using the given
uncertainty ellipsoid at the $k$th step. Since the attitude dynamics
is nonlinear, the admissible boundary of the state at the $(k+l)$th
step is not an ellipsoid in general. We assume that the uncertainty
ellipsoid at the $k$th step is sufficiently small that states in the
uncertainty ellipsoid can be approximated using the linearized equations
of motion. 

 \textit{Center:} For the given center $(\hat{R}_k,\hat{\Omega}_k)$,
the center of the uncertainty ellipsoid at step $(k+l)$ is
$(\hat{R}_{k+l}^{f},\hat{\Omega}_{k+l}^{f})$ obtained using the
discrete equations of motion, \refeqn{findf0}, \refeqn{updateR0},
and \refeqn{updatew0}:
\begin{gather}
h S(J\hat{\Omega}_k+\frac{h}{2} \hat{M}_k) = \hat{F}_k J_d - J_d
\hat{F}_k^T,\label{eqn:findf}\\
\hat{R}_{k+1}^{f} = \hat{R}_k \hat{F}_k,\label{eqn:updateR}\\
J\hat{\Omega}_{k+1}^{f} = \hat{F_k}^T
\hat{\Omega}_k+\frac{h}{2}\hat{F_k}^T \hat{M_k} +\frac{h}{2}
\hat{M}_{k+1}.\label{eqn:updatePi}
\end{gather}
This integrator yields a map $(\hat R_k,\hat\Omega_k)\mapsto(\hat
R^f_{k+1},\hat\Omega_{k+1}^f)$, and this process is repeatedly applied
to find the center at the $(k+l)$th step, $(\hat R^f_{k+l},
\hat\Omega_{k+l}^f)$.

\textit{Uncertainty matrix:}
At the $(k+1)$th step, the state is represented by
perturbations from the center $(\hat R^f_{k+1},
\hat\Omega_{k+1}^f)$:
\begin{align*}
    R_{k+1}&=\hat{R}_{k+1}^{f} e^{S(\zeta_{k+1}^{f})},\\
    \Omega_{k+1}&=\hat{\Omega}_{k+1}^{f}+\delta\Omega_{k+1}^{f},
\end{align*}
for some $\zeta_{k+1}^{f},\delta\Omega_{k+1}^{f}\in\Re^3$. The
uncertainty matrix at the $(k+1)$th step is obtained by finding a
bound on $\zeta_{k+1}^{f},\delta\Omega_{k+1}^{f}\in\Re^3$. Assume
that the uncertainty ellipsoid at the $k$th step is sufficiently
small. Then, $\zeta_{k+1}^{f},\delta\Omega_{k+1}^{f}$ are
represented by the following linear equations in~\cite{pro:acc06}
\begin{align*}
x_{k+1}^{f} & = A_k^f x_k,
\end{align*}
where $x_k=[\zeta_k;\delta\Omega_k]\in\Re^6$, and
$A_k^f\in\Re^{6\times 6}$ can be suitably defined. Since
$(R_k,\Omega_k)\in\mathcal{E}(\hat{R}_k,\hat{\Omega}_k,P_k)$,
$x_k\in\mathcal{E}_{\Re^6}(0,P_k)$ by the definition of the
uncertainty ellipsoid given in \refeqn{ueso}. Then we can show that
$A_k^f x_k$ lies in
\begin{align*}
A_k^f x_k&\in\mathcal{E}_{\Re^6}\!\parenth{0,A_k^f P_k
\parenth{A_k^f}^T}.
\end{align*}
Thus, the uncertainty matrix at the $(k+1)$th step is given by
\begin{align}
P_{k+1}^f & = A_k^f P_k \parenth{A_k^f}^T.\label{eqn:Pkpf}
\end{align}
In summary, the uncertainty ellipsoid at the $(k+l)$th step is
computed using \refeqn{findf}, \refeqn{updateR}, \refeqn{updatePi},
and \refeqn{Pkpf} as:
\begin{align}\label{eqn:flow}
    (R_{k+l},\Omega_{k+l})\in\mathcal{E}(\hat{R}_{k+l}^{f},
\hat{\Omega}_{k+l}^{f},P_{k+l}^f).
\end{align}

\subsection{Measurement update}
The measurement update finds an uncertainty ellipsoid in the
state space using the measurements and sensor error models. The
measured attitude and the angular velocity define the center of the
measurement uncertainty ellipsoid, and the sensor error models give
the uncertainty matrix.

\textit{Center:} The center of the uncertainty ellipsoid,
$(\hat{R}_{k+l}^{m},\hat\Omega_{k+l}^m)$ is obtained from
measurements. Let the measured directions to the known points be
$\tilde B_{k+l}=[\tilde b^1,\cdots,\tilde b^m]\in\Re^{3\times m}$.
Then, the attitude $\hat{R}_{k+l}^{m}$ satisfies the following
necessary condition given in \refeqn{san2006}
\begin{gather}
\parenth{\hat{R}_{k+l}^m}^T\tilde L_{k+l}-\tilde L_{k+l}^T
\hat{R}_{k+l}^m=0,\label{eqn:meaR}
\end{gather}
where $\tilde L_{k+l}=E_{k+l} W_{k+l}\tilde{B}_{k+l}^T\in\Re^{3\times 3}$.
The attitude matrix 
is given by a QR factorization of $\tilde L_{k+l}$ as in \refeqn{Rhat}
\begin{align}
\hat{R}_{k+l}^m=\parenth{Q_q\sqrt{(Q_rQ_r^T)^{-1}}Q_q^T} \tilde
L_{k+l},\label{eqn:meaRex}
\end{align}
where $Q_q\in\SO$ is an orthogonal matrix and $Q_r\in\Re^{3\times
3}$ is a upper triangular matrix satisfying $\tilde L_{k+l}=Q_qQ_r$.

The angular velocity is measured directly by
\begin{align}
\hat{\Omega}_{k+l}^{m}= \tilde\Omega_{k+l}.\label{eqn:meaw}
\end{align}

\textit{Uncertainty matrix:} We represent the actual state at the
$(k+l)$th step as perturbations from the measured center:
\begin{align}
    R_{k+l}&=\hat{R}_{k+l}^{m} e^{S(\zeta_{k+l}^{m})},\label{eqn:Rkpm}\\
    \Omega_{k+l}&=\hat{\Omega}_{k+l}^{m}+\delta\Omega_{k+l}^{m},\label{eqn:Omegakpm}
\end{align}
for $\zeta_{k+l}^{m},\delta\Omega_{k+l}^{m}\in\Re^3$. The uncertainty
matrix is obtained by finding a bound on $\zeta_{k+l}^{m},
\delta\Omega_{k+l}^{m}$.

We transform the uncertainties in measuring the body directions to
known fixed points into uncertainties in the rotation matrix by
\refeqn{meaR}. Using the error model in \refeqn{bi}, the actual
directions corresponding to ${B}_{k+l}$ are given by
\begin{align}
{B}_{k+l} & = \tilde{B}_{k+l} +
\delta\tilde{B}_{k+l},\label{eqn:Bkl}
\end{align}
where $\delta B_{k+l}=\bracket{S(\nu^1)\tilde
b^1,\cdots,S(\nu^m)\tilde b^m}\in\Re^{3\times m}$.

The actual directions ${B}_{k+l}$ and the actual attitude $R_{k+l}$
at the $(k+l)$th step also satisfy \refeqn{meaRex};
\begin{gather}
R_{k+l}^T L_{k+l}- L_{k+l}^T{R}_{k+l}=0,\label{eqn:meaRac}
\end{gather}
where $L_{k+l}=E_{k+l} W_{k+l}{B}_{k+l}^T\in\Re^{3\times 3}$.
Substitute \refeqn{Rkpm} and \refeqn{Bkl} into \refeqn{meaRac},
and use 
$S(x)A+A^TS(x)=S(\braces{\tr{A}I_{3\times 3}-A}x)$ for
$A\in\Re^{3\times 3},x\in\Re^3$, to get: 
\begin{align*}
\braces{\tr {\parenth{\hat{R}_{k+l}^m}^T
\tilde{L}_{k+l}}-\parenth{\hat{R}_{k+l}^m}^T \tilde{L}_{k+l}}
\zeta_{k+l}^m \qquad\qquad\quad \\=-\sum_{i=1}^m w_i
\braces{\tr{\tilde{b}^i (e^i)^T \hat{R}_{k+l}^m}I_{3\times
3}-\tilde{b}^i (e^i)^T \hat{R}_{k+l}^m}\nu^i.
\end{align*}
We can rewrite the above equation as
\begin{align}\label{eqn:zetakpm}
\zeta_{k+l}^m &= \sum_{i=1}^m \mathcal{A}_{k+l}^{m,i}\nu^i,
\end{align}
where $\mathcal{A}_{k+l}^{m,i}\in\Re^{3\times 3}$ is defined
appropriately.

The perturbation of the angular velocity $\delta\Omega_{k+l}^{m}$ is
equal to the angular velocity measurement error $\upsilon_{k+l}$,
\begin{align}\label{eqn:delPikpm}
\delta\Omega_{k+l}^m = \upsilon_{k+l}.
\end{align}

Define
$x_{k+l}^m=\bracket{\zeta^m_{k+l};\,\delta\Omega^m_{k+l}}\in\Re^6$.
Using \refeqn{zetakpm} and \refeqn{delPikpm},
\begin{align*}
x_{k+l}^m & = H_1 \sum_{i=1}^m \mathcal{A}_{k+l}^{m,i} \nu^i_{k+l} +
H_2 \upsilon_{k+l},
\end{align*}
where $H_1=[I_{3\times 3},\, 0_{3\times 3}]^T,H_2=[0_{3\times 3},\,
I_{3\times 3}]^T\in\Re^{6\times 3}$. This expresses $x_{k+l}^m$ as
a linear combination of the sensor errors $\nu^i$ and $\upsilon$.
Using the measurement uncertainties \refeqn{Sk} and \refeqn{Tk},
we can show that the terms in the right hand side of the above
equation are in the following uncertainty ellipsoids:
\begin{align*}
H_1 \mathcal{A}_{k+l}^{m,i} \nu^i_{k+l} & \in
\mathcal{E}_{\Re^6}\parenth{0,H_1 \mathcal{A}_{k+l}^{m,i}S^i_{k+l}
\parenth{\mathcal{A}_{k+l}^{m,i}}^TH_1^T},\\
H_2 \upsilon_{k+l} & \in \mathcal{E}_{\Re^6}\parenth{0,H_2 T_{k+l}
 H_2^T}.
\end{align*}
Thus, the uncertainty ellipsoid for $x_{k+l}^{m}$ is obtained as the
vector sum of the above uncertainty ellipsoids. The measurement
update obtains a minimal ellipsoid that contains the vector sum of
these uncertainty ellipsoids. Using expressions for such a minimal
ellipsoid 
given in~\cite{jo:MaNo1996} and~\cite{jo:DuWaPo2001}, we get: 
\begin{align}
P_{k+l}^m & = \braces{\sum_{i=1}^m \sqrt{\tr{P^{m,i}_{k+l,R}}}
+\sqrt{\tr{P^m_{k+l,\Omega}}}}\nonumber\\
& \quad \times \braces{\sum_{i=1}^m
\frac{P^{m,i}_{k+l,R}}{\sqrt{\tr{P^{m,i}_{k+l,R}}}}
+\frac{P^m_{k+l,\Omega}}{\sqrt{\tr{P^m_{k+l,\Omega}}}}},\label{eqn:Pmkl}
\end{align}
where
\begin{align*}
P^{m,i}_{k+l,R}&=H_1
\mathcal{A}_{k+l}^{m,i}S^i_{k+l}\parenth{\mathcal{A}_{k+l}^{m,i}}^TH_1^T,\\
P^m_{k+l,\Omega}&=H_2 T_{k+l} H_2^T.
\end{align*}
In summary, the measured uncertainty ellipsoid at the $(k+l)$th step
is defined by \refeqn{meaRex}, \refeqn{meaw}, and \refeqn{Pmkl};
\begin{align}\label{eqn:mea}
    (R_{k+l},\Omega_{k+l})\in\mathcal{E}(\hat{R}_{k+l}^{m},
\hat{\Omega}_{k+l}^{m},P_{k+l}^m).
\end{align}

\subsection{Filtering procedure}
The filtering procedure obtains a new uncertainty ellipsoid
compatible with both the predicted and the measured uncertainty
ellipsoids. From \refeqn{flow} and \refeqn{mea}, we know that: 
\begin{align*}
    (R_{k+l},\Omega_{k+l})\in\mathcal{E}(\hat{R}_{k+l}^{f},
\hat{\Omega}_{k+l}^{f},P_{k+l}^f)\bigcap
    \mathcal{E}(\hat{R}_{k+l}^{m},\hat{\Omega}_{k+l}^{m},P_{k+l}^m).
\end{align*}
The intersection of two ellipsoids is not generally an ellipsoid,
and it is inefficient to describe an irregular subset in the
multidimensional space numerically. We find a minimal uncertainty
ellipsoid containing this intersection.
We omit the subscript $(k+l)$ here for convenience.

The measurement uncertainty ellipsoid, $\mathcal{E}(\hat{R}^{m},
\hat{\Omega}^{m},P^m)$, is identified by its center $(\hat{R}^{m},
\hat{\Omega}^{m})$, and the uncertainty ellipsoid in $\Re^6$:
\begin{align}
(\zeta^m,\delta\Omega^m)\in\mathcal{E}_{\Re^6}(0_{6\times
1},P^m),\label{eqn:ellRm}
\end{align}
where $S(\zeta^m)=\mathrm{logm} \parenth{(\hat{R}^m)^T R}\in\so$,
$\delta\Omega^m=\Omega-\hat{\Omega}^m\in\Re^3$. Similarly, the
flow uncertainty ellipsoid, $\mathcal{E}(\hat{R}^{f},\hat{\Omega}^{f},
P^f)$, is identified by its center $(\hat{R}^{f},\hat{\Omega}^{f})$,
and the uncertainty ellipsoid in $\Re^6$:
\begin{align}
(\zeta^f,\delta\Omega^f)\in\mathcal{E}_{\Re^6}(0_{6\times
1},P^f),\label{eqn:ellRf}
\end{align}
where $S(\zeta^f)=\mathrm{logm} \parenth{(\hat{R}^f)^T R}\in\so$,
$\delta\Omega^f=\Omega-\hat{\Omega}^f\in\Re^3$. An element
$(R^f,\Omega^f)\in\mathcal{E}(\hat{R}^{f},\hat{\Omega}^{f},P^f)$ is
given by
\begin{align}
R^f & = \hat{R}^{f} e^{S(\zeta^f)},\label{eqn:Rf}\\
\Omega^f & = \hat{\Omega}^{f} + \delta{\Omega}^f.\label{eqn:Pif}
\end{align}
Define $\hat\zeta^{mf},\delta\hat\Omega^{mf}\in\Re^3$ such that
\begin{align}
\hat{R}^f&=\hat{R}^m e^{S(\hat\zeta^{mf})},\label{eqn:errz}\\
\hat{\Omega}^f& =\hat{\Omega}^m
+\delta\hat\Omega^{mf}.\label{eqn:errPi}
\end{align}
Thus, $\hat\zeta^{mf},\delta\hat\Omega^{mf}$ represent the
difference between the centers of the two ellipsoids.

Substituting \refeqn{errz}, \refeqn{errPi} into \refeqn{Rf},
\refeqn{Pif}, we obtain
\begin{align}
R^f & = \hat{R}^m e^{S(\hat\zeta^{mf})} e^{S(\zeta^f)},\nonumber\\
& \simeq \hat{R}^m e^{S(\hat\zeta^{mf}+\zeta^f)},\\
\Omega^f & = \hat{\Omega}^m + \parenth{\delta\hat\Omega^{mf} +
\delta{\Omega}^f},
\end{align}
where we assumed that $\hat\zeta^{mf}, \zeta^f$ are sufficiently
small. Thus, the uncertainty ellipsoid obtained by the flow update,
$\mathcal{E}(\hat{R}^{f},\hat{\Omega}^{f},P^f)$ is identified by the
measured $(\hat{R}^m,\hat{\Omega}^m)$ and the following uncertainty
ellipsoid in $\Re^6$:
\begin{align}
\mathcal{E}_{\Re^6}( \hat{x}^{mf} ,P^f),\label{eqn:ellRmf}
\end{align}
where
$\hat{x}^{mf}=\bracket{\hat\zeta^{mf};\delta\hat\Omega^{mf}}\in\Re^6$.

We seek a minimal ellipsoid that contains the intersection: 
\begin{align}
\mathcal{E}_{\Re^6}(0_{6\times 1},P^m)\bigcap
\mathcal{E}_{\Re^6}(\hat{x}^{mf}
,P^f)\subset\mathcal{E}_{\Re^6}(\hat{x},P),
\end{align}
where $\hat{x}=[\hat\zeta;\delta\hat\Omega]\in\Re^6$. 
Using the expression for a minimal ellipsoid containing the
intersection of two ellipsoids presented in~\cite{jo:MaNo1996},
$\hat{x}$ and $P$ are given by
\begin{align*}
\hat{x}&=L\hat{x}^{mf},\\
P&=\beta(q) (I-L)P^m,
\end{align*}
where
\begin{align*}
\beta(q) & = 1 + q - (\hat{x}^{mf})^T (P^{m})^{-1} L \hat{x}^{mf},\\
L & = P^{m} (P^{m} + q^{-1} P^f)^{-1}.
\end{align*}
The constant $q$ is chosen to minimize $\mathrm{tr}\!\bracket{P}$.
We convert $\hat{x}$ to points in $\T\SO$ using the common center
$(\hat{R}^{m},\hat{\Omega}^{m})$.

In summary, the uncertainty ellipsoid at $(k+l)$th step is
\begin{align}
    (R_{k+l},\Omega_{k+l})\in\mathcal{E}(\hat{R}_{k+l},\hat{\Omega}_{k+l},P_{k+l}),
\end{align}
where
\begin{equation}
\hat{R}_{k+l}=\hat{R}_{k+l}^m e^{S(\hat\zeta)}, \;\
\hat{\Omega}_{k+l} = \hat{\Omega}_{k+l}^m + \delta\hat\Omega,\;\
P_{k+l} = P.
\end{equation}

\subsection{Properties of the estimator}

The steps outlined above are repeated to get a dynamic filter.
This attitude estimator has
no singularities since the attitude is represented by a rotation
matrix. Orthogonality of the rotation matrix is preserved as it is
updated by the structure-preserving Lie group variational integrator.
This estimator can be used for highly nonlinear large angle
maneuvers of a rigid body. It is also robust to the distribution of
the sensor noise since we only use ellipsoidal bounds on the noise.
The measurements need not be periodic, the estimation is repeated
whenever new measurements become available. 
We can also extend this attitude estimator to the case when
angular velocity measurements are not available. The filtering step
is modified to find an intersection of the non-degenerate predicted
uncertainty ellipsoid and the degenerate measurement uncertainty
ellipsoid.

\section{Numerical Simulation}
Numerical simulation results are presented for estimation of the
attitude dynamics of an uncontrolled rigid spacecraft in a circular
orbit about a large central body, including gravity gradient effects.
The on orbit spacecraft model is given in~\cite{pro:acc06}.

The mass, length and time dimensions are normalized by the
spacecraft mass, the maximum length of the spacecraft, and the
orbital angular velocity, respectively. The inertia of the
spacecraft is chosen as ${J}=\mathrm{diag}\bracket{1,\,2.8,\,2}$.
The maneuver is an arbitrary large attitude change completed in a
quarter of the orbit. The initial conditions are chosen as
\begin{alignat*}{2}
R_0 & = \mathrm{diag}[-1,-1,1],&
\Omega_0&=[2.316,\,0.446,\,-0.591]\,\mathrm{rad/s},\\
\hat R_0 & = I_{3\times 3},&\quad
\hat{\Omega}_0&=[2.116,\,0.546,\,-0.891]\,\mathrm{rad/s}.
\end{alignat*}
The corresponding initial estimation errors are
$\norm{\zeta_0}=180\,\mathrm{deg}$,
$\norm{\delta\Omega_0}=21.43\frac{\pi}{180}\,\mathrm{rad/s}$. Note
that the actual initial attitude is opposite to the estimated
initial attitude. The initial uncertainty matrix is given by
\begin{align*}
P_0 = 2\,\mathrm{diag}\bracket{\pi^2[1,\,1,\,1],\,
\parenth{\frac{\pi}{6}}^2[1,\,1,\,1]},
\end{align*}
so that $x_0^TP_0^{-1}x_0=0.7553\leq 1$.

\begin{figure}
    \centerline{\subfigure[Estimation error $\norm{\zeta_k}$, $\norm{\delta\Omega_k}$]{
    \includegraphics[width=0.50\columnwidth]{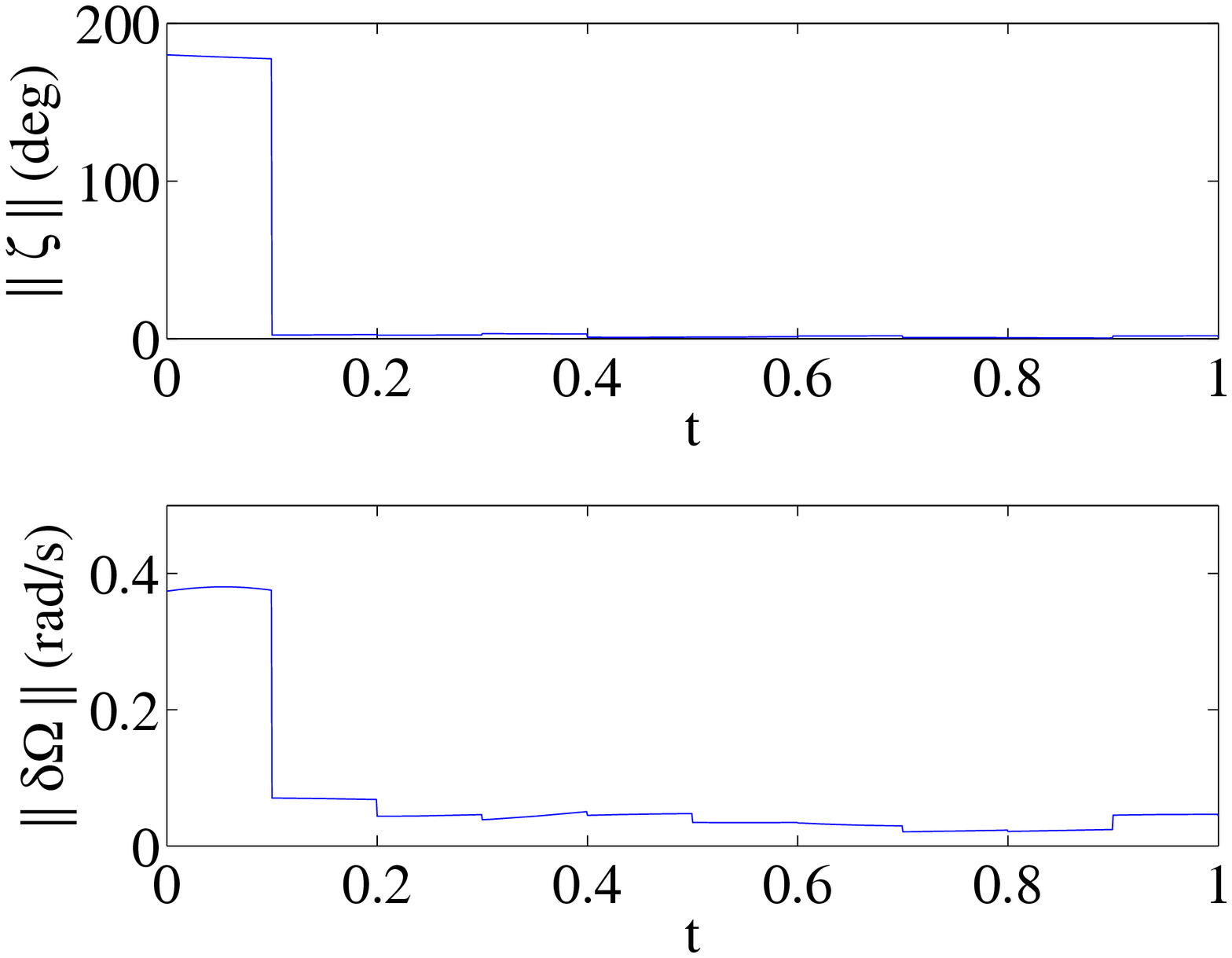}}
    \hfill
    \subfigure[Size of uncertainty $\tr{P_k}$]{
    \includegraphics[width=0.48\columnwidth]{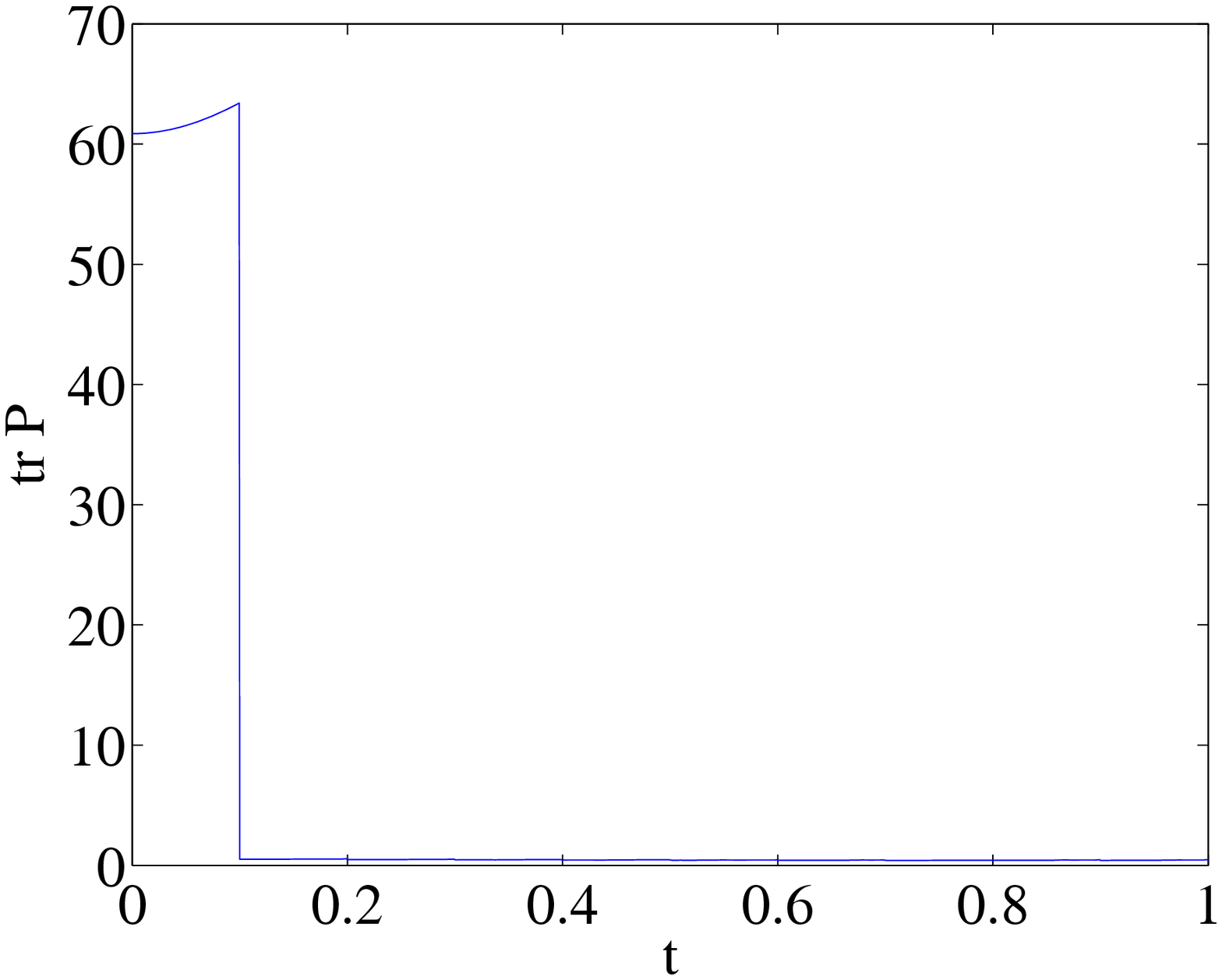}}
    }
    \caption{Estimation with attitude and angular velocity
    measurement}\label{fig:full}
    \centerline{\subfigure[Estimation error $\norm{\zeta_k}$, $\norm{\delta\Omega_k}$]{
    \includegraphics[width=0.50\columnwidth]{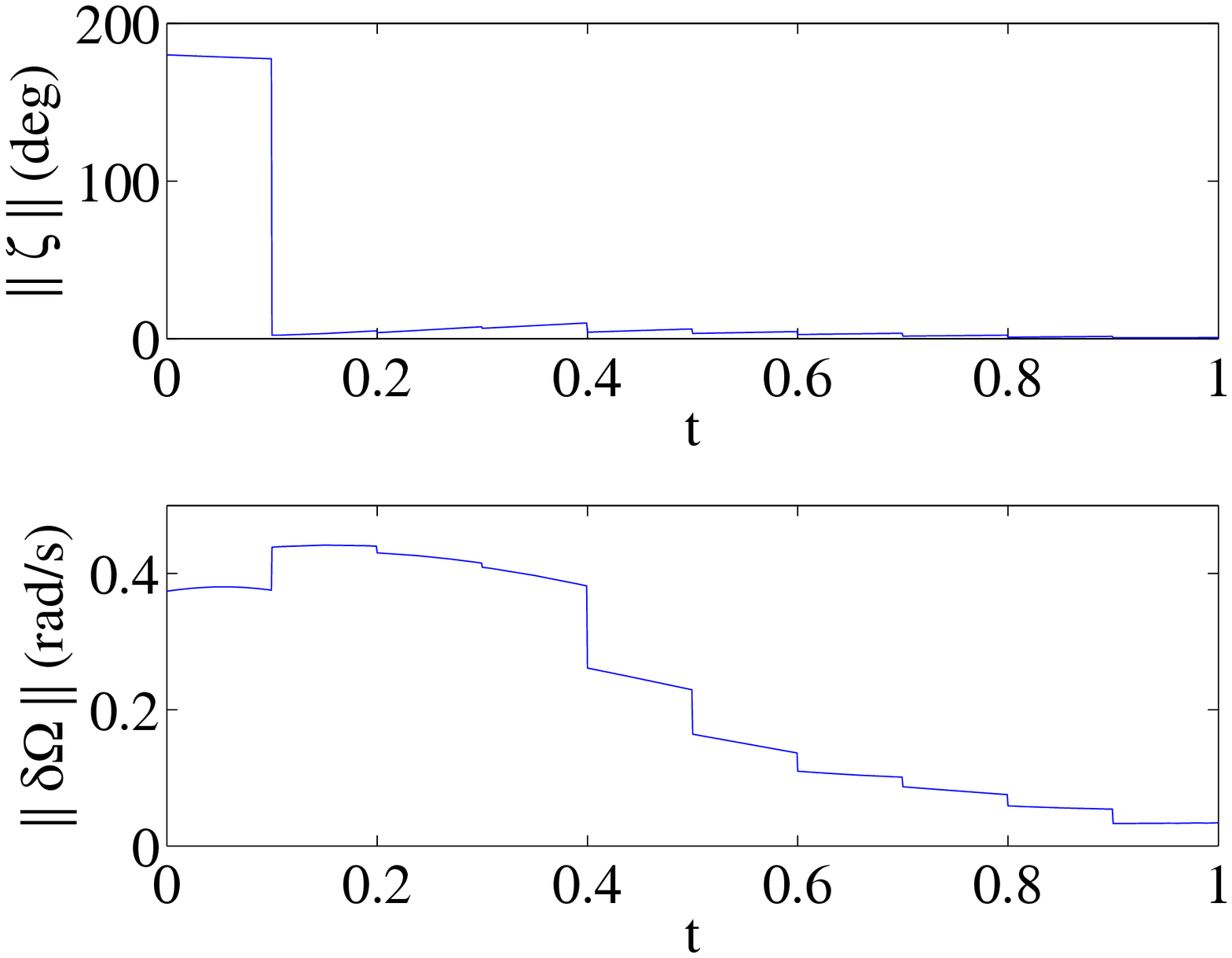}}
    \hfill
    \subfigure[Size of uncertainty $\tr{P_k}$]{
    \includegraphics[width=0.48\columnwidth]{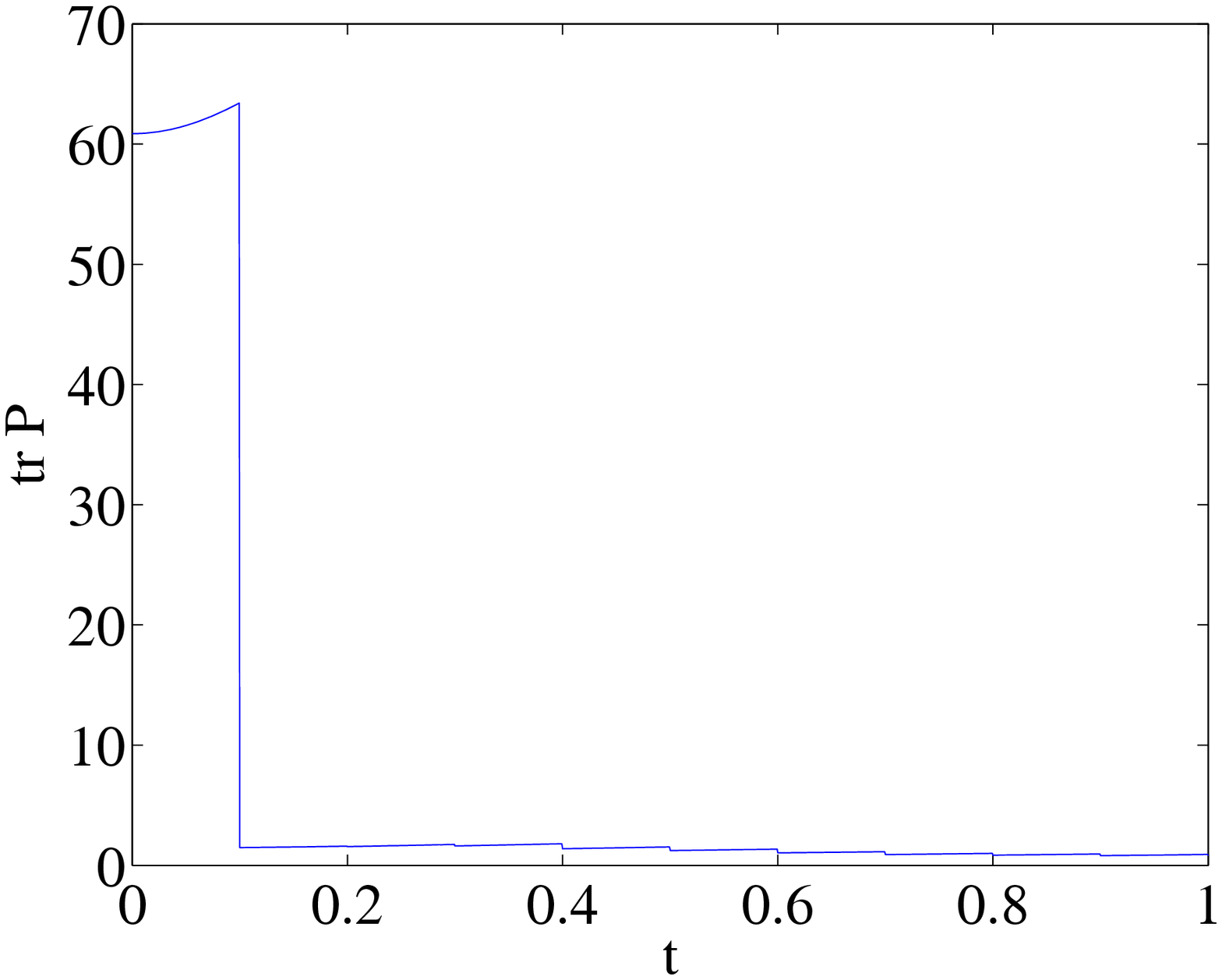}}
    }
    \caption{Estimation with attitude measurement only}\label{fig:att}
\end{figure}

We assume that the measurements are available ten times per quarter
orbit. The measurement noise is assumed to be normally distributed with
uncertainty matrices given by
\begin{align*}
S^i_k=\parenth{7\frac{\pi}{180}}^2 I_{3\times
3}\,\mathrm{rad^2},\quad T_k=\parenth{7\frac{\pi}{180}}^2 I_{3\times
3}\,\mathrm{rad^2/s^2}.
\end{align*}
We consider two cases. \reffig{full} shows simulation results when
both the attitude and the angular velocity are measured. \reffig{att}
shows simulation results when angular velocity measurements are not
available. In each figure, the left plot shows the attitude and
angular velocity estimation errors, and the right plot shows the size
of the uncertainty ellipsoid. The estimation errors and the size of
uncertainty decrease rapidly after the first measurement. When
the angular velocity measurements are not available, the estimation
error for the angular velocity converges relatively slowly as seen
in \reffig{att}.(a). For both cases, the terminal attitude error,
and the terminal angular velocity error are less than
$0.88\,\mathrm{deg}$, and $0.04\,\mathrm{rad/s}$, respectively.

\section{Conclusion}
A deterministic estimator for the attitude dynamics of a rigid body
in a potential field with bounded measurement errors is presented.
An uncertainty ellipsoid is obtained at each estimation step, and
the dynamics is propagated using Lie group variational integrators.
The center of the uncertainty ellipsoid is the point estimate, and
its size determines the accuracy of the estimate. The estimation
scheme is optimal in the sense that the size of the uncertainty is
minimized at each estimation step. It is also global and robust to
the distribution of measurement noise. This estimator can be
extended to include the effects of process noise and to the case
when only attitude measurements are available.
These extensions are not described in this paper.

\bibliography{ue}

\begin{thebibliography}{10}
\providecommand{\url}[1]{#1}
\csname url@rmstyle\endcsname
\providecommand{\newblock}{\relax}
\providecommand{\bibinfo}[2]{#2}
\providecommand\BIBentrySTDinterwordspacing{\spaceskip=0pt\relax}
\providecommand\BIBentryALTinterwordstretchfactor{4}
\providecommand\BIBentryALTinterwordspacing{\spaceskip=\fontdimen2\font plus
\BIBentryALTinterwordstretchfactor\fontdimen3\font minus
  \fontdimen4\font\relax}
\providecommand\BIBforeignlanguage[2]{{%
\expandafter\ifx\csname l@#1\endcsname\relax
\typeout{** WARNING: IEEEtran.bst: No hyphenation pattern has been}%
\typeout{** loaded for the language `#1'. Using the pattern for}%
\typeout{** the default language instead.}%
\else
\language=\csname l@#1\endcsname
\fi
#2}}

\bibitem{jo:wahba}
G.~Wahba, ``A least squares estimate of satellite attitude, {P}roblem 65-1,''
  \emph{SIAM Review}, vol.~7, no.~5, p. 409, 1965.

\bibitem{BaOs.AES85}
I.~Y. Bar-Itzhack and Y.~Oshman, ``Attitude determination from vector
  observations; quaternion estimation,'' \emph{IEEE Transactions on Aerospace
  and Electronic Systems}, vol.~21, no.~1, pp. 128--136, 1985.

\bibitem{Shu.JAS90}
M.~D. Shuster, ``Kalman filtering of spacecraft attitude and the
  {Q}{U}{E}{S}{T} model,'' \emph{Journal of the Astronautical Sciences},
  vol.~38, no.~3, pp. 377--393, 1990.

\bibitem{ReHu.Au04}
H.~Rehbinder and X.~Hu, ``Drift-free attitude estimation for accelerated rigid
  bodies,'' \emph{Automatica}, vol.~40, no.~4, pp. 653--659, 2004.

\bibitem{Pi.JGCD00}
M.~L. Psiaki, ``Attitude determination filtering via extended quaternion
  estimation,'' \emph{AIAA Journal of Guidance, Control and Dynamics}, vol.~23,
  no.~2, pp. 206--214, 2000.

\bibitem{jo:TheSkaSou.IEEETSP94}
Y.~Theodor, U.~Shaked, and C.~E. de~Souza, ``A game theory approach to robust
  discrete-time {$H_\infty$}-estimation,'' \emph{IEEE Transactions on Signal
  Processing}, vol.~42, no.~6, pp. 1486--1495, 1994.

\bibitem{jo:Sc1968}
F.~C. Schweppe, ``Recursive state estimation: Unknown but bounded errors and
  system inputs,'' \emph{IEEE Transactions on Automatic Control}, vol.~13,
  no.~1, pp. 22--28, 1968.

\bibitem{jo:MaNo1996}
D.~G. Maksarov and J.~P. Norton, ``State bounding with ellipsoidal set
  description of the uncertainty,'' \emph{International Journal of Control},
  vol.~65, no.~5, pp. 847--866, 1996.

\bibitem{jo:DuWaPo2001}
C.~Durieu, E.~Walter, and B.~Polyak, ``Multi-input multi-output ellipsoidal
  state bounding,'' \emph{Journal of Optimization Theory and Applications},
  vol. 111, no.~2, pp. 273--303, 2001.

\bibitem{jo:solwahba}
J.~L. Farrell, J.~C. Stuelpnagel, R.~H. Wessner, J.~R. Velman, and J.~E. Brock,
  ``A least squares estimate of satellite attitude, {S}olution 65-1,''
  \emph{SIAM Review}, vol.~8, no.~3, pp. 384--386, 1966.

\bibitem{jo:ShOh1981}
M.~D. Shuster and S.~D. Oh, ``Three-axis attitude determination from vector
  observations,'' \emph{Journal of Guidance Control and Dynamics}, vol.~4,
  no.~1, pp. 70--77, 1981.

\bibitem{pro:san2006}
A.~K. Sanyal, ``Optimal attitude estimation and filtering without using local
  coordinates, {P}art {I}: {U}ncontrolled and deterministic attitude
  dynamics,'' in \emph{Proceedings of the American Control Conference}, 2006,
  pp. 5734--5739.

\bibitem{pro:acc06}
T.~Lee, M.~Leok, and N.~H. McClamroch, ``Attitude maneuvers of a rigid
  spacecraft in a circular orbit,'' in \emph{Proceedings of the American
  Control Conference}, 2006, pp. 1742--1747.

\bibitem{pro:cca05}
------, ``A {L}ie group variational integrator for the attitude dynamics of a
  rigid body with applications to the 3{D} pendulum,'' in \emph{Proceedings of
  the IEEE Conference on Control Applications}, 2005, pp. 962--967.

\end{thebibliography}
\bibliographystyle{IEEEtran}

\end{document}